# Enhanced Active Power Filter Control for Nonlinear Non-Stationary Reactive Power Compensation

Phen Chiak See, Vin Cent Tai, Marta Molinas, Kjetil Uhlen, and Olav Bjarte Fosso

*Norwegian University of Science and Technology*

June 18, 2012

**Abstract.** This paper describes a method to implement Reactive Power Compensation (RPC) in power systems that possess nonlinear non-stationary current disturbances. The Empirical Mode Decomposition (EMD) introduced in the Hilbert-Huang Transform (HHT) is used to separate the disturbances from the original current waveform. These disturbances are subsequently removed. Following that, Power Factor Correction (PFC) based on the well-known $p$-$q$ power theory is conducted to remove the reactive power. Both operations were implemented in a shunt Active Power Filter (APF). The EMD significantly simplifies the singulation and the removal of the current disturbances. This helps to effectively identify the fundamental current waveform. Hence, it simplifies the implementation of RPC on nonlinear non-stationary power systems.

## 1   Introduction

Currently, several methods have been proposed to characterize the nonlinear non-stationary power disturbances in the power systems. For instance, several methods that identify non-stationary disturbances were mentioned in [Lobos and Rezmer, 1997; Janik et al., 2007; Lobos et al., 2007; Lobos et al., 2008; and Su et al., 2009] (to name a few). They are based mainly on the wavelet analysis and Prony's estimation. Among them, [Janik et al., 2007; Lobos et al., 2007; Lobos et al., 2008; and Su et al., 2009] have studied the interesting transient switching behavior of Power Factor Correction (PFC) capacitors in wind turbines, as well as its effect on the performance of power systems. Although these methods excel in detecting the non-stationary power disturbances, they do not suggest the practical way to remove them prior to the implementation of PFC.

A sample current waveform with nonlinear non-stationary disturbances from the wind turbine (as calculated by [Lobos et al., 2007]) is shown in Figure 1. It was mentioned in [Lobos et al., 2007] that the existence of such disturbances complicates the implementation of Reactive Power Compensation (RPC) with capacitor banks. In the authors' opinion, it is important to create an adaptive RPC strategy that could handle the random onset of such disturbances. Furthermore, the method should be able to deal with the non-existence of zero mean in the distorted current waveforms. A reliable method that could separate the distorted current into a fundamental waveform and several distinguished disturbances is needed.

Therefore, the authors propose to use the Empirical Mode Decomposition (EMD) method (which was developed as part of the Hilbert-Huang Transform (HHT) [Huang et al., 1998]) to effectively identify and singulate the nonlinear non-stationary current disturbances from the power system. An instantaneous RPC method that combines the EMD and the shunt Active Power Filter (APF) control based on $p$-$q$ power theory [Akagi et al., 1983; Akagi et al. 1984; and Watanabe et al., 2004] was developed. The study further validates the feasibility of this strategy on an unbalanced 3-phase 4-wire power system model that contains distorted currents.

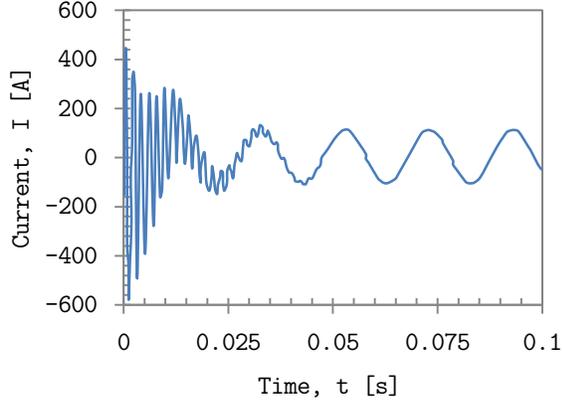

Figure 1. Calculated current caused by capacitor bank switching at 8m/s wind speed (adapted from [Lobos et al., 2007]).

This paper begins with an introduction on the concept of the proposed method. Following that, descriptions on the steps involved in EMD that singulates the current disturbances are given. Subsequently, the EMD-enhanced control strategy for shunt APF developed based on the $p$-$q$ power theory is addressed. The results obtained from the proposed method are compared with the RPC results obtained from the existing shunt APF control so as to validate its performance. The paper finally ends with conclusions.

## 2  Empirical mode decomposition

The EMD was introduced by [Huang et al., 1998] as an innovative method to decompose the nonlinear non-stationary data into several instantaneous oscillating components known as the Intrinsic Mode Functions (IMFs) that possess zero mean. It enables the use of the Hilbert Transform (HT) to analyze the energy-time-frequency characteristics of the data.

Generally, EMD involves a finite number of sifting processes. Each of them begins with the identification of the mean between the upper envelope and the lower envelope. As shown in Figure 2, the upper envelope is an interpolated cubic spline that connects all local maxima in the time series data. Similarly, the lower envelope is an interpolated cubic spline that connects the local minima. The mean between the upper envelope and the lower envelope is designated as $m_i$. Next, the $m_i$ are subtracted from the time series data to obtain the $i$-th component of sift, $h_i$:

$$X(t) - m_i = h_i \quad (1)$$

The sifting process is then repeated on $h_i$ until the component of sift exhibits the following characteristics: (i.) possesses zero mean, and (ii.) doesn't possess either positive local minima values or negative local maxima.

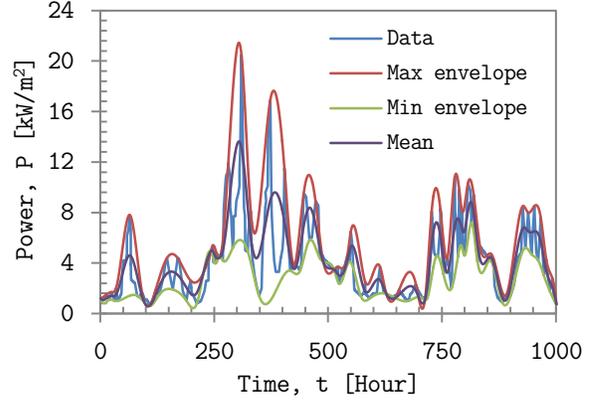

Figure 2. Sample showing how the maximum envelope, minimum envelope and the mean values of time series data are defined in EMD.

The sift component is then designated as the $k$-th component of IMF, $c_k$. The $c_k$ is subtracted from the original time series data to obtain the $k$-th residue, $r_k$:

$$X(t) - c_k = r_k \quad (2)$$

The Standard Deviation (SD) is used as the stopping criteria for the $i$-th level's iteration, which can be fixed between 0.2 and 0.3 [Huang et al., 1998]:

$$SD = \sum_{t=1}^{T}\left[\frac{|h_{i-1}(t) - h_i(t)|^2}{h_{i-1}^2(t)}\right] \quad (3)$$

From the authors' experience, this value determines the intensity of IMF subtraction in (2). Larger portion from the original data will be removed as IMFs if the SD value is high.

The sifting process at $i$-th and $k$-th levels are repeated until a pre-determined stopping criteria is reached at step $n$ (i.e., when the $r_n$ contain least (or no) instantaneous frequency oscillation). As a result, we have obtained a number of IMFs ($c_0, c_1, … c_n$) that represent the instantaneous oscillating components of the data and a final residue, $r_n$ decomposed from the original time series data:

$$X(t) = \sum_{k=1}^{n} c_k + r_n \quad (4)$$

These components respectively represent the fundamental waveform, as well as the instantaneous oscillating components of the nonlinear non-stationary data. This allows us to use the HT for interpreting the time series data as follows:

$$x(t) = \Re\left\{\sum_{j=1}^{n} a_j(t) e^{i \int \omega_j(t) dt}\right\} \quad (5)$$

with, $a_j(t)$ denotes the magnitude of HT's analytical function, $Z(t) = X(t) + iY(t)$, and $\omega = \frac{d\theta(t)}{dt}$. $Y(t)$ and $\theta(t)$ are HT value and the phase of $Z(t)$.

In the RPC developed in this work, $r_n$ is used as the fundamental waveform as it contains the least instantaneously oscillating component. The IMFs represent the localized nonlinear non-stationary disturbances that will be removed during the compensation process.

## 3 Proposed reactive power compensation

The authors propose to combine EMD with the instantaneous RPC based on $p$-$q$ theory for effective PFC. A sample test case with unbalanced load was created by the authors as shown in Figure 3.

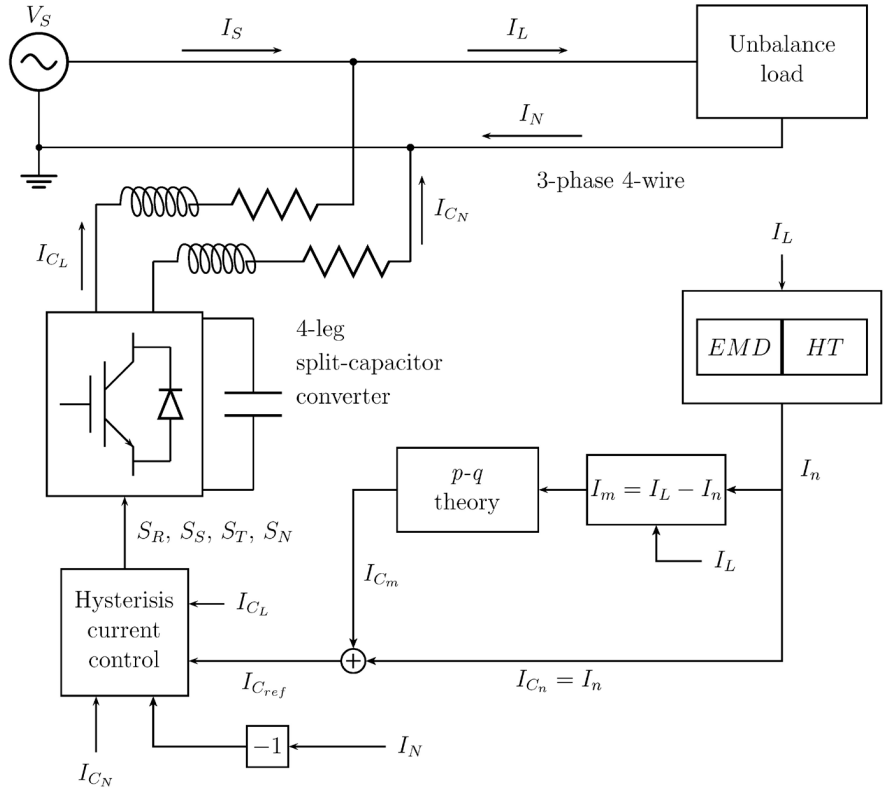

Figure 3. Sample test case of the proposed RPC method.

## 3.1 Power disturbances in the test case

As shown in Figure 4, the load current at line 1 of the simulation model contains non-stationary disturbances (from 0.088s to 0.094s as shown in Figure 5). Meanwhile, nonlinear current distortions are observed in load currents at line 3. The voltages obtained from the model are shown in Figure 6. The distorted currents and the line voltages create severe distortions in instantaneous power (as shown in Figure 7).

The authors argue that the current distortion at line 1 is similar with the current distortion shown in Figure 1 (e.g., from time 0.02s to 0.04s). The nature of distortions in both cases is same although they have different time characteristics. The exactly similar distortion was illustrated in Figure 8 of [Gaing, 2003].

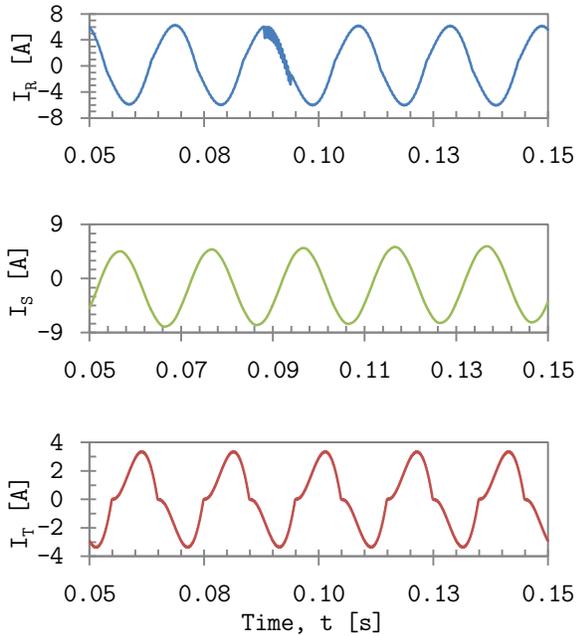

Figure 4. From top to bottom – load current at (a.) line 1; (b.) line 2; and (c.) line 3, within the simulation model.

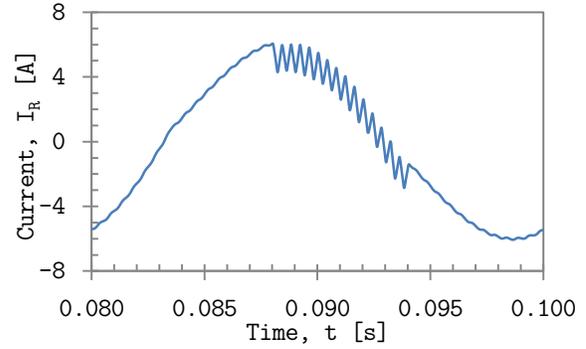

Figure 5. Load current disturbance detected at line 1 (from 0.08s to 0.1s).

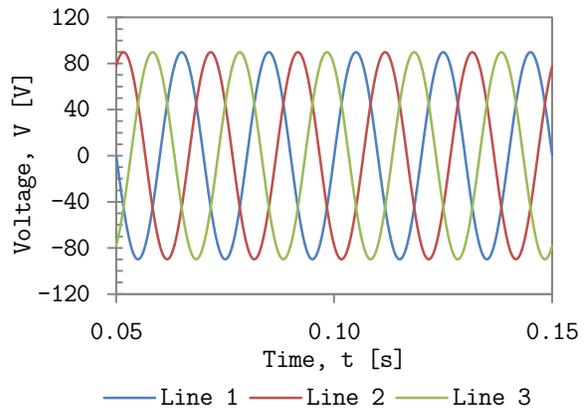

Figure 6. Load voltages at all lines within the simulation model.

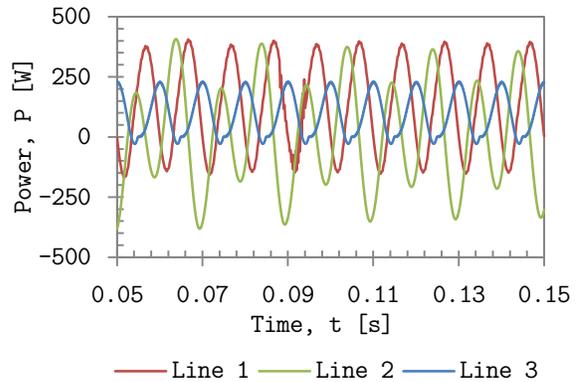

Figure 7. Instantaneous power at all lines within the simulation model.

The nonlinear current distortion at line 3 resembles the ordinary current distortion found in the unbalanced nonlinear power systems.

## 3.2 Detection of current disturbances

The EMD decomposes the distorted currents (at line 1 and line 3, respectively) to: (i.) $r_n$ that represents the fundamental current waveform, and (ii.) multiple IMFs that represent the current disturbances. The $r_n$ and IMFs are treated as distinguished components in the shunt APF, which simultaneously implements the disturbance removal and PFC.

As shown in Figure 8b, two disturbances in current at line 1 are identified. These components are removed from the original current waveform to yield the fundamental current waveform (Figure 8c). Similarly, two disturbances in current at line 3 are identified (see Figure 9b for the current disturbances and Figure 10 for detailed information on the IMF 1 shown in Figure 9b). The fundamental current waveform is shown in Figure 9c. These disturbances are removed using a 4-leg split-capacitor converter controlled by a hysteresis current controller in the shunt APF.

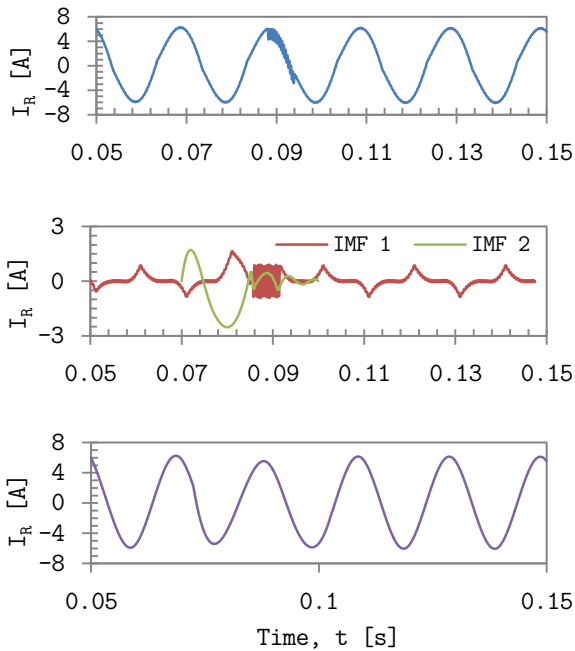

Figure 8. From top to bottom – (a.) load current at line 1; (b.) current disturbances at line 1 detected as IMFs by EMD; and (c.) smoothed current at line 1 after the removal of disturbances.

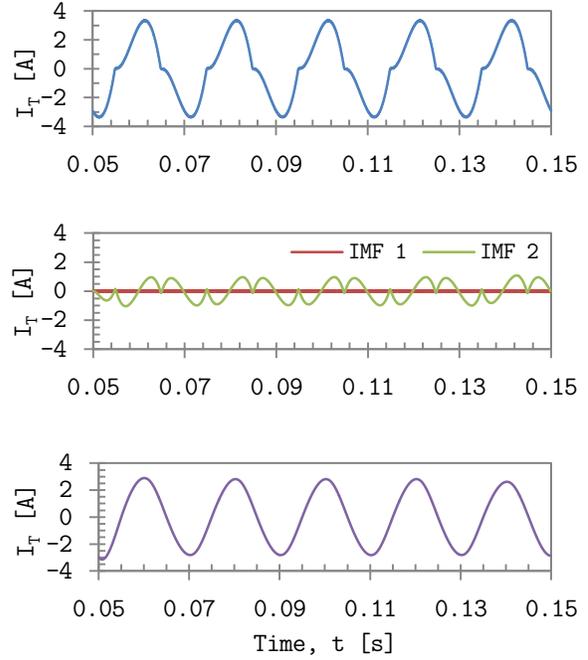

Figure 9. From top to bottom – (a.) load current at line 3; (b.) current disturbances at line 3 detected as IMFs by EMD; and (c.) smoothed current at line 3 after the removal of disturbances.

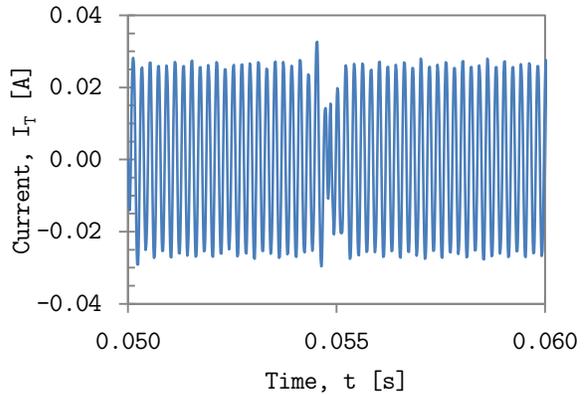

Figure 10. Current disturbances detected at line 3 (IMF 1).

## 3.3 Control of shunt APF

The proposed EMD-enhanced shunt APF control was derived from hybridizing the EMD with the existing instantaneous RPC based on the modified *p-q* theory described in [Akagi et al., 1999; and Togasawa et al., 1994]. The latter is chosen because of its capability to eliminate the zero sequence current and zero sequence power that exist in the 3-phase 4-wire systems under unbalanced loading condition. The main objective of this hybridization is to determine the compensating currents that have to be produced by the shunt APF in the presence of nonlinear non-stationary currents.

The RPC starts with the measurement of load currents (generalized as the line current, $I_L$ (see Figure 3), which consists of the $I_R$, $I_S$, $I_T$ components) and the neutral current, $I_N$. The measured $I_L$ are used as the input for EMD. The current disturbances, $I_n$ detected by EMD, as well as the fundamental current waveforms, $I_m$ gained from subtracting the current disturbances from the respective $I_L$ are then used to control of the shunt APF.

In the control model, the compensating currents for eliminating the $I_n$ are injected back to the power system as the compensating currents, $I_{C_n}$ for eliminating the disturbances. Afterwards, the modified $p$-$q$ theory is used to calculate the compensating currents, $I_{C_m}$ for PFC. This calculation is rather straightforward as the $I_m$ doesn't contain any nonlinear non-stationary disturbances. The $I_{C_m}$ and $I_{C_n}$ are then summed to obtain the reference compensation currents, $I_{C_{ref}}$.

Next, a hysteresis current controller is used to generate the control signals ($S_R$, $S_S$, $S_T$, and $S_N$) for the 4-leg split-capacitor converter based on the $I_{C_{ref}}$. The converter injects the compensating load currents, $I_{C_L}$ back to the power systems to achieve the complete RPC.

On the other hand, the neutral current, $I_N$ is compensated by simply feeding the neutral line with the neutral current in the reverse direction. At the end of the current injections, the results anticipated are distortion-free sinusoidal source currents with unity power factor and zero reactive power, $q$ from the power source.

## 4 Simulation results

The proposed RPC with EMD-enhanced shunt APF control (henceforth known as the EMD-enhanced RPC) was implemented on the simulation model. The results obtained are shown in Figure 11.

In order to compare and validate its performance, RPC using the shunt APF control method derived from the existing modified $p$-$q$ theory (henceforth known as the RPC without EMD) is as well provided in the same figure. The steps involved in implementing the RPC without EMD is similar to [Akagi et al., 1999; and Togasawa et al. 1994].

As shown, the EMD-enhanced RPC demonstrates several advantages over the RPC without EMD:

- **Power output:** both methods are able to produce zero reactive power, $q$. However, the EMD-enhanced RPC creates stable active power, $p$ and $q$ as compared to the RPC without EMD. Non-stationary oscillation in $p$ and $q$ (between 0.075s and 0.1s) are not observed in the EMD-enhanced RPC.
- **Power factor:** both methods yield unity power factor. Nevertheless, the oscillation in power factor is observed in the RPC without EMD due to the existence of non-stationary current components. The EMD-enhanced RPC doesn't possess this problem.
- **Total Harmonic Distortion (THD) index:** both methods achieve convergence in the THD index. This indicates that the harmonic distortion in both cases is improving over time. There is no significant difference/improvement created by the EMD-enhanced RPC in this respect.
- **Current waveform:** the nonlinear current distortions are successfully removed by both RPC methods. Nevertheless, the RPC without EMD failed to remove the non-stationary current disturbance. In fact, this current disturbance causes the non-stationary oscillations observed in the $p$, the $q$ and the power factor. As for the EMD-enhanced shunt APF control, the non-stationary current disturbance is perfectly removed.

In short, the EMD-enhanced shunt APF control excels in compensating 3-phase 4-wire power systems with nonlinear non-stationary current disturbances. In particular, the advantage of using it on the non-stationary current disturbances is very apparent.

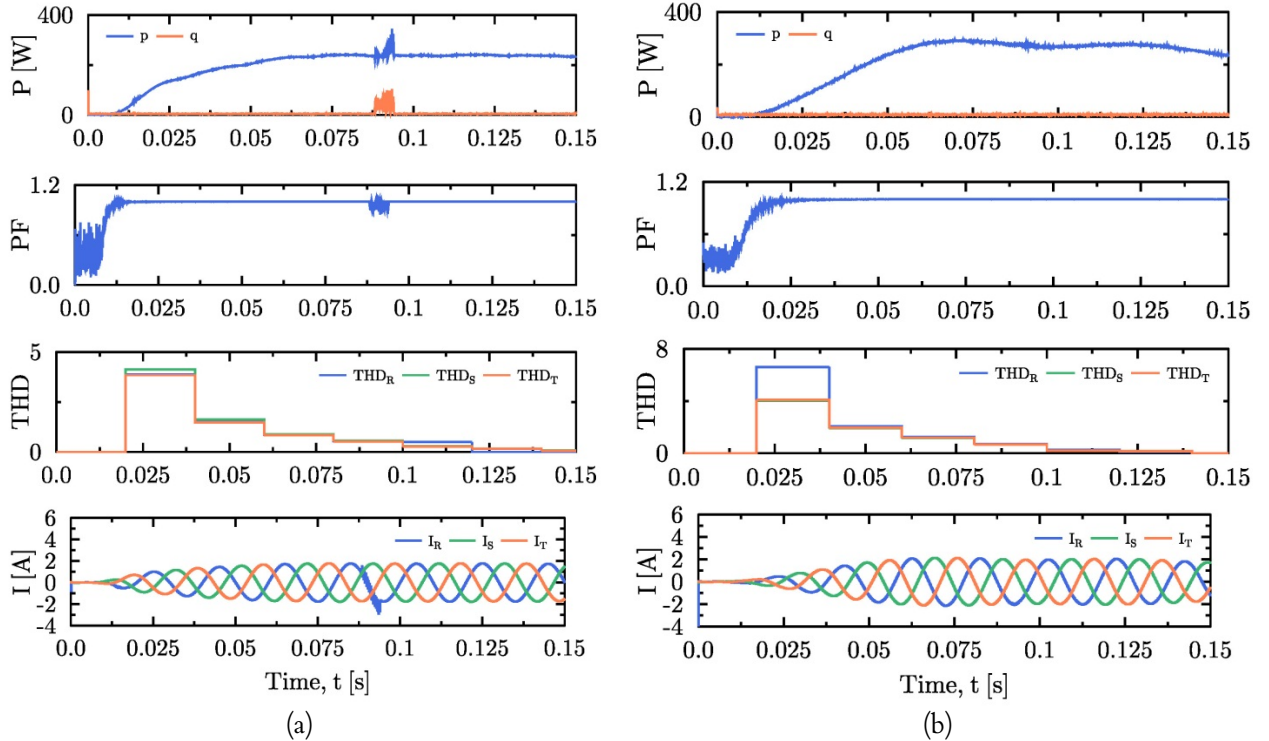

Figure 11. Simulation results obtained from: (a.) the existing shunt APF control, and (b.) the EMD-enhanced shunt APF control.

## 5  Conclusions

This paper has presented a shunt APF control method for RPC in 3-phase 4-wire power systems that contain nonlinear non-stationary power disturbances. A hybrid control method that combines the EMD and the existing shunt APF control based on the modified *p-q* theory was developed to serve the purpose. Steps involved in implementing the proposed RPC are addressed and assessment is performed to validate its performance. The results have proven that the proposed method excels in removing non-stationary current disturbances, while performing equally well with the existing shunt APF control method when handling the nonlinear current disturbances. Note-worthily, the EMD is able to produce control input for singulation of current disturbances in real time. The performance of this feature will be examined at the next level of study.

In future, this model can be extended to solve other power system problems such as the control of active damping used in the Wide Area Monitoring (WAM) of smart-grids. This is because, the inter-area oscillating component found in the WAM problems are nonlinear non-stationary in nature. The use of the control method based on Fourier Analysis (FA) is no longer valid in such situation due to spectral leakage problems. Also, it is unable to precisely identify the energy-time-frequency information in time domain.

## 6  Acknowledgements

The authors wish to thank the Norwegian University of Science and Technology for all supports given in this research.

## References

Akagi H., Kanazawa Y. and Nabae A. Generalized theory of the instantaneous reactive power in three-phase circuits. In Proceedings of the International Power Electronic Conference, pp. 1375-1386. Tokyo, Japan, March 27-31, 1983.

Akagi H., Kanazawa Y. and Nabae A. Instantaneous reactive power compensator comprising switching devices without energy storage components. IEEE Transactions on Industry Applications 20 (3): 625-630, 1984.

Akagi H., Satoshi O. and Kim H. The theory of instantaneous power in three-phase four-wire systems: a comprehensive approach. In Conference Record of the 1999 IEEE Industry Applications Conference, pp. 431-439. Arizona, USA, October 3-7, 1999.

Gaing Z. L. Implementation of power disturbance classifier using wavelet-based neural networks. In Proceedings of the IEEE Bologna PowerTech Conference. Bologna, Italy, June 23-26, 2003.

Huang N. E., Shen Z., Long S. R., Wu M. C., Shih H. H., Zheng Q., Yen N. C., Tung C. C. and Liu H. H. The empirical mode decomposition and the Hilbert spectrum for nonlinear and non-stationary time series analysis. Proceedings of the Royal Society 454: 903-995, 1998.

Janik P., Lobos T., Rezmer J., Waclawek Z. and Thiringer T. Wind generator transients' computation using Prony method. In Proceedings of the International Conference on Clean Electrical Power, pp. 598-604. Capri, Italy, May 21-23, 2007.

Lobos T. and Rezmer J. Real-time determination of power system frequency. IEEE Transactions on Instrumentation and Measurement 46 (4): 877-881, 1997.

Lobos T., Rezmer J., Janik P. and Waclawek Z. Prony and nonlinear regression methods used for determination of transient parameters in wind energy conversion system. In Proceedings of the PowerTech Conference, pp. 1764-1769. Lausanne, Switzerland, July 1-5, 2007.

Lobos T., Rezmer J., Sikorski T. and Waclawek Z. Power distortion issues in wind turbine power systems under transient states. Turkish Journal of Electrical Engineering 16 (3): 229-238, 2008.

Su S. P., Luo X. and Qin Z. Q. Study on transient poewr quality detection of grid-connected wind power generation system based on wavelet transform. In Proceedings of the International Conference on Energy and Environmental Technology, pp. 861-864. Guilin, China, October 16-18, 2009.

Togasawa S., Murase T., Nakanoand H. and Nabae A. Reactive power compensation based on a novel cross-vector theory. Transactions on Industrial Application of the IEEE Japan 114 (3): 340-341, 1994.

Watanabe E. H., Aredes M. and Akagi H. The *p-q* theory for active filter control: some problems and solutions. Revista Controle and Automacao 15 (1): 78-84, 2004.